\documentclass[11pt]{article}
\usepackage{amscd}
\usepackage{amsfonts}
\usepackage{amsmath}
\usepackage{amssymb}
\usepackage{amsthm}
\usepackage{bbm}
\usepackage{CJK}
\usepackage{fancyhdr}
\usepackage{graphicx}
\usepackage{hyperref}
\usepackage{indentfirst}
\usepackage{latexsym}
\usepackage{mathrsfs}
\usepackage{xypic}

\newtheorem{theorem}{Theorem}[section]
\newtheorem{lemma}[theorem]{Lemma}
\newtheorem{definition}[theorem]{Definition}
\newtheorem{proposition}[theorem]{Proposition}
\newtheorem{example}[theorem]{Example}

\newtheorem{corollary}[theorem]{Corollary}

\usepackage[top=1in,bottom=1in,left=1.25in,right=1.25in]{geometry}
\def\<{\langle}
\def\>{\rangle}
\def\a{\alpha}
\def\b{\beta}

\def\c{\cdot}

\def\g{\gamma}

\def\o{\otimes}

\date{}
\begin{document}
\renewcommand{\baselinestretch}{1.2}
\renewcommand{\arraystretch}{1.0}
\title{\bf On split regular  BiHom-Poisson superalgebras}
\author{{\bf Shuangjian Guo$^{1}$, Yuanyuan Ke$^{2}$\thanks{Corresponding author, Email: keyy086@126.com } }\\
{\small 1. School of Mathematics and Statistics, Guizhou University of Finance and Economics} \\
{\small  Guiyang  550025, P. R.  China} \\
{\small 2. School of Mathematics and Computer Science, Jianghan University}\\
 {\small   Wuhan, 430056, P. R. China}}
 \maketitle
 \maketitle
\begin{center}
\begin{minipage}{13.cm}

{\bf \begin{center} ABSTRACT \end{center}}
The paper introduces the class of split regular BiHom-Poisson superalgebras, which is a natural
generalization of split regular Hom-Poisson algebras and  split regular BiHom-Lie superalgebras.  By developing techniques of connections
of roots for this kind of algebras, we show that such a split regular BiHom-Poisson superalgebras $A$ is of the form  $A=U+\sum_{\a}I_\a$ with $U$ a subspace of a maximal abelian
subalgebra $H$ and any $I_{\a}$, a well described ideal of $A$, satisfying $[I_\a, I_\b]+I_\a I_\b = 0$ if
$[\a]\neq [\b]$. Under  certain conditions, in the case of  $A$ being of maximal length, the
simplicity of the algebra is characterized.


 \end{minipage}
 \end{center}
 \normalsize\vskip1cm

\section*{INTRODUCTION}
\def\theequation{0. \arabic{equation}}
\setcounter{equation} {0}

The notion of Hom-Lie algebras was first introduced by Hartwig,
Larsson and Silvestrov in \cite{Hartwig}, who developed an approach to deformations 
of the Witt and  Virasoro algebras basing on $\sigma$-deformations.
In fact, Hom-Lie algebras include Lie algebras as a subclass, but the deformation  of Lie algebras twisted by a homomorphism.

The twisting  parts of the defining identities can transfer one algebra to the other algebraic structures.
In \cite{Makhlouf2008,Makhlouf2010},  Makhlouf and  Silvestrov introduced the notions of
 Hom-associative algebras, Hom-coassociative coalgebras, Hom-bialgebras  and Hom-Hopf algebras.
The original definition of a Hom-bialgebra involved two linear maps, one was twisting the
associativity condition and the other was twisting the coassociativity condition. In the case of Hom-Lie algebras, the relevant structure for a tensor theory is a Hom-Poisson algebra structure. A Hom-Poisson algebra has simultaneously a Hom-Lie algebra structure and a Hom-associative algebra structure, satisfying the Hom-Leibniz identity in  \cite{Makhlouf10}.  In \cite{Wang15}, Wang, Zhang and Wei characterized Hom-Leibniz superalgebras and Hom-Leibniz Poisson superalgebras, and presented the methods to construct these superalgebras.

A BiHom-algebra is an algebra in such a way that the identities defining the structure
are twisted by two homomorphisms $\phi$ and $\psi$. This class of algebras was introduced from a
categorical approach in \cite{Graziani} which as an extension of the class of Hom-algebras. If the two
linear maps are the same automorphisms, BiHom-algebras will be return to Hom-algebras.
These algebraic structures include BiHom-associative algebras, BiHom-Lie algebras and
BiHom-bialgebras. The representation theory of BiHom-Lie algebras was introduced by
Cheng and Qi in \cite{Cheng16}, in which, BiHom-cochain complexes, derivation, central extension, derivation
extension, trivial representation and adjoint representation of BiHom-Lie algebras were
studied. More applications of  BiHom-algebras, BiHom-Lie superalgebras,  BiHom-Lie colour algebras and BiHom-Novikov algebras
can be found in (\cite{Liu17}, \cite{Wang16}, \cite{Abdaoui17}, \cite{Guo2018}).

The class of the split algebras is specially related to addition quantum numbers, graded contractions and deformations. For instance, for a physical system
which displays a symmetry, it is interesting to know the detailed structure of the split
decomposition, since its roots can be seen as certain eigenvalues which are the additive
quantum numbers characterizing the state of such system. Determining the structure of
split algebras will become more and more meaningful in the area of research in mathematical physics. Recently, in (\cite{Albuquerque 2018}-\cite{Calderon2016}, \cite{Cao2015}-\cite{Cao17},  \cite{Zhang2018}-\cite{Zhang2017}), the structure of different classes of split algebras have been determined by the techniques of connections of roots. The purpose
of this paper is to consider the class of  split regular BiHom-Poisson superalgebras, which is a natural extension of   split regular BiHom-Lie superalgebras
and   split Hom-Lie  superalgebras.

In Section 2, we prove that such an arbitrary split regular BiHom-Poisson superalgebras $A$ is of the form  $A=U+\sum_{\a}I_\a$ with $U$ a subspace of a maximal abelian subalgebra $H$ and any $I_{\a}$, a well described ideal of $A$, satisfying $[I_\a, I_\b]+I_\a I_\b = 0$ if
$[\a]\neq [\b]$.

In Section 3, we present that under certain conditions, in the case of $A$ being of maximal
length, the simplicity of the algebra is characterized.

\section{Preliminaries}
\def\theequation{\arabic{section}.\arabic{equation}}
\setcounter{equation} {0}

Throughout this paper, we will denote by $\mathbb{N}$ the set of all nonnegative integers and by $\mathbb{Z}$ the set of all integers.  Split regular BiHom-Poisson superalgebras  are considered of arbitrary dimension and over an arbitrary base field  $\mathbb{K}$.
 Any unexplained definitions and notations can be found in \cite{Graziani} and \cite{Makhlouf10}, and  we recall some basic definitions and results related to our paper.

\noindent{\bf 1.1. Hom-Poisson superalgebra }
A Hom-Poisson superalgebra is a Hom-Lie superalgebra  $(A,[\cdot, \cdot],\phi)$ endowed with a Hom-associative superproduct,that is, a bilinear product denoted by juxtaposition such that
\begin{eqnarray*}
\phi(x) (yz)=(xy)\phi(y),
\end{eqnarray*}
for all $x,y,z\in A$, and such that the Hom-Leibniz superidentity
\begin{eqnarray*}
[xy,\phi(z)]=\phi(x)[y,z]+ (-1)^{\overline{j}\overline{k}}[x,z]\phi(y)
\end{eqnarray*}
holds for any $x\in A_{\overline{i}}, y\in A_{\overline{j}},z\in A_{\overline{k}}$ and $\overline{i},\overline{j},\overline{k}\in \mathbb{Z}_2$.

If $\phi$ is furthermore a Poisson automorphism, that is,  a linear bijective on such that $\phi([x,y])=[\phi(x),\phi(y)]$ and $\phi(xy)=\phi(x)\phi(y)$ for any $x,y\in A$, then $A$ is called a regular Hom-Poisson superalgebra.

\noindent{\bf 1.2. BiHom-associative superalgebra} A  BiHom-associative superalgebra is a 4-tuple $(A,\mu,\phi,\psi)$,
where $A$ is a  superspace, $\phi: A\rightarrow A$, $\psi: A\rightarrow A$ and
 $\mu: A\o A\rightarrow A$  are linear maps, with notation $\mu(a\o b)=ab$,
 satisfying the following conditions, for all  $a,a',a''\in A$:
\begin{eqnarray*}
&&\phi\circ\psi=\psi\circ\phi,\\
&&\phi(aa')=\phi(a)\phi(a'),\psi(aa')=\psi(a)\psi(a'),\\
&&\phi(a)(a'a'')=(aa')\psi(a'').
\end{eqnarray*}
And the maps $\phi$ and $\psi$ are called the structure maps of $A$.
\smallskip

Clearly, a Hom-associative algebra $(A,\mu,\phi)$ can be regarded as the BiHom-associative
algebra $(A,\mu,\phi,\phi)$.
\medskip

\noindent{\bf 1.3. BiHom-Lie superalgebra }  A BiHom-Lie superalgebra $L$ is a $\mathbb{Z}_2$-graded algebra $L=L_{\overline{0}}\oplus L_{\overline{1}}$, endowed with an even bilinear mapping $[\c,\c]: L\times L\rightarrow L$ and two homomorphisms $\phi, \psi: L\rightarrow L$
 satisfying the following conditions, for all  $x\in L_{\overline{i}},y\in L_{\overline{j}},z\in L_{\overline{k}}$ and $\overline{i},\overline{j},\overline{k}\in \mathbb{Z}_2$:
\begin{eqnarray*}
&& [x,y]\subset L_{\overline{i}+\overline{j}},\\
&&\phi\circ\psi=\psi\circ\phi,\\
&&[\psi(x),\phi(y)]=-(-1)^{\overline{i}\overline{j}}[\psi(y),\phi(x)],\\
&&(-1)^{\overline{k}\overline{i}}[\psi^{2}(x),[\psi(y),\phi(z)]]+(-1)^{\overline{i}\overline{j}}[\psi^{2}(y),[\psi(z),\phi(x)]]+(-1)^{\overline{j}\overline{k}}[\psi^{2}(z),[\psi(x),\phi(y)]]=0.
\end{eqnarray*}

When $\phi$ and $\psi$ are algebra automorphisms, it is said  that $L$ is a regular BiHom-Lie superalgebra.

\section{Decomposition}
\def\theequation{\arabic{section}. \arabic{equation}}
\setcounter{equation} {0}

\begin{definition}
A BiHom-Poisson superalgebra $A$ is a BiHom-Lie superalgebra  $(A,[\cdot, \cdot],\phi,\psi)$ endowed with a BiHom-associative superproduct, that is, a bilinear product denoted by juxtaposition such that
\begin{eqnarray*}
\phi(x) (yz)=(xy)\psi(y),
\end{eqnarray*}
for all $x,y,z\in A$, and such that the BiHom-Leibniz superidentity
\begin{eqnarray*}
[xy,\phi\psi(z)]=\phi(x)[y,\phi(z)]+(-1)^{|j||k|}[x,\psi(z)]\phi(y)
\end{eqnarray*}
holds for any $x\in A_{\overline{i}}, y\in A_{\overline{j}},z\in A_{\overline{k}}$ and $\overline{i},\overline{j},\overline{k}\in \mathbb{Z}_2$.
\end{definition}
Furthermore, if $\phi$ and $\psi$ are Poisson automorphisms, it is said that $A$ is  a regular BiHom-Poisson superalgebra.

 \begin{example}
 Let $(A, \mu, [\cdot, \cdot])$ be a Poisson superalgebra and $\phi, \psi: A\rightarrow A$ two commuting Possion superalgebras automorphism. If we endow the underlying linear space $A$ with  new products $[\cdot,\cdot]':A\times A\rightarrow A, \mu'=\mu\circ (\phi\o \psi)$ defined by $[x,y]'=[\phi(x),\psi(y)], \mu'(x\o y)=\phi(x)\psi(y)$ for any $x,y\in A$, we know that $(A, \mu',  [\cdot, \cdot]',  \phi, \psi)$ becomes a regular BiHom-Poisson superalgebra.
\end{example}

\begin{example}
Let $A=A_{\overline{0}}\oplus A_{\overline{1}}$ be a 2-dimensional superspace, where $A_{\overline{0}}$ is generated by $e_1$ and $A_{\overline{1}}$ is generated by $e_2$ and nonzero products and $\phi, \psi: A\rightarrow A$ are given by
\begin{eqnarray*}
\phi(e_1)=ae_1~~\phi(e_2)=ae_2~~\psi(e_1)=be_1~~\psi(e_2)=be_2,~~ a, b\neq 0,1,\\
e_1\c e_1=abe_1, ~~e_2\c e_2= abe_1,~~ e_1\c e_2=e_2\c e_1=ab e_2, ~~[e_2, e_2]=2ab e_2.
\end{eqnarray*}
Then $(A, \c, [\c,\c], \phi, \psi)$ is a regular BiHom-Poisson superalgebra.

\end{example}

\begin{example}
Let $A=A_{\overline{0}}\oplus A_{\overline{1}}$ be a 3-dimensional superspace, where $A_{\overline{0}}$ is generated by $e_1, e_2$ and $A_{\overline{1}}$ is generated by $e_3$ and nonzero products  and $\phi, \psi: A\rightarrow A$  are given by
\begin{eqnarray*}
\phi(e_1)=ae_1~~\phi(e_2)=e_1+e_2~~\phi(e_3)=0, ~~a\neq 0,1,\\
\psi(e_1)=be_1~~\psi(e_2)=e_1+e_2, ~~\psi(e_3)=0, ~~b\neq 0,1,\\
e_1\c e_2=ae_1, ~~e_2\c e_2= e_1+e_2,~~[e_1, e_2]=ae_1.
\end{eqnarray*}
Then $(A, \c, [\c,\c], \phi, \psi)$ is a regular BiHom-Poisson superalgebra.
\end{example}

Note that $A_{\overline{0}}$ is a BiHom-Poisson algebra called the even or bosonic part of $A$, while $A_{\overline{1}}$ is called the odd or fermonic part of $A$. The
usual regularity concepts will be understood in the graded sense. That is, a subalgebra $H=H_{\overline{0}}\oplus H_{\overline{1}}$ of $A$ is a graded  subspace such that $[H,H]+HH\subset A$ and $\phi(H)=\psi(H)=H$. A graded subspace $I=I_{\overline{0}}\oplus I_{\overline{1}}$ of $A$ is called an ideal if $[I,A]+IA+AI\subset I$ and $\phi(I)=\psi(I)=I$. A BiHom-Poisson superalgebra $A$ will be called simple if $[A, A]+AA\neq 0$  and its only ideals are $\{0\}$ and $A$.

We recall from \cite{Zhang2018} that a BiHom-Lie superalgebra $(A, [\cdot, \cdot], \phi, \psi)$ and a maximal abelian sualgebra $H$ of $A$, for a linear functional
\begin{eqnarray*}
\a:H_{\overline{0}}\rightarrow \mathbb{K},
\end{eqnarray*}
we define the root space of $A$ associated to $\a$ as the subspace
\begin{eqnarray*}
A_{\a}:=\{v_{\a}\in A:[h_{\overline{0}},\phi(v_{\a})]=\a(h)\phi\psi(v_{\a}), \mbox{for any $h_{\overline{0}}\in H_{\overline{0}}$}\}.
\end{eqnarray*}
The elements $\a:H_{\overline{0}}\rightarrow \mathbb{K}$ satisfying $A_{\a}\neq 0$ are called roots of $A$ with respect to $H$ and we denote $\Lambda:=\{\a\in (H_{\overline{0}})^{\ast}/\{0\}:A_{\a}\neq 0\}$. We call that $A$ is a split regular BiHom-Lie superalgebra with respect to $H$ if
\begin{eqnarray*}
A=H\oplus \bigoplus_{\a\in \Lambda}A_{\a}.
\end{eqnarray*}
We also say that $\Lambda$ is the root system of $A$.

To simplify notation, the mappings $\phi|_H, \psi|_H, \phi^{-1}|_H, \psi^{-1}|_H: H\rightarrow H$ will be denoted by $\phi, \phi^{-1}$, $\psi$ and $\psi^{-1}$ respectively.

We recall   some properties of split regular BiHom-Lie superalgebras  that can be found in \cite{Zhang2018}.
\begin{lemma}
 Let $(A, [\cdot, \cdot], \phi, \psi)$ be a split regular BiHom-Lie superalgebra. Then,
for any $\a,\b\in A\cup \{0\}$,

(1)$\phi(A_{\a})=A_{\a\phi^{-1}}$ and $\psi(A_{\a})=A_{\a\psi^{-1}}$,

(2) $\phi^{-1}(A_{\a})=A_{\a\phi}$ and $\psi^{-1}(A_{\a})=A_{\a\psi}$,

(3) $[A_{\a},A_{\b}]\subset A_{\a\phi^{-1}+\b\psi^{-1}}$,

(4) If $\a\in \Lambda$, then $\a \phi^{-z_1}\psi^{-z_2}\in \Lambda$ for any $z_1, z_2 \in \mathbb{Z}$,

(5) $A_{0}=H$.
\end{lemma}

\begin{lemma}
 Let $A$ be a split regular BiHom-Poisson superalgebra. Then
for any $\a,\b\in A\cup \{0\}$, we have $A_{\a}A_{\b}\subset A_{\a\phi^{-1}+\b\psi^{-1}}$.
\end{lemma}

{\bf Proof} Let $h_{\overline{0}}\in H_{\overline{0}}, v_\a\in A_{\a, {\overline{i}}}$ and $v_\b\in A_{\b, {\overline{j}}}$, we can write
\begin{eqnarray*}
[h_{\overline{0}}, \phi(v_\a v_\b)]=[\psi\psi^{-1}(h_{\overline{0}}),\phi(v_\a v_\b)],
\end{eqnarray*}
and denote $h'_{\overline{0}}=\psi^{-1}(h_{\overline{0}})$. By applying the BiHom-Leibniz superidentity,  we get
\begin{eqnarray*}
[\psi\psi^{-1}(h_{\overline{0}}),\phi(v_\a v_\b)]&=&[\psi(h'_{\overline{0}}),\phi(v_\a v_\b)]\\
&=&-[\psi(v_\a v_\b),\phi(h'_{\overline{0}})]\\
&=& -(-1)^{\overline{0}\overline{j}}[\psi(v_\a), h']\phi\psi(v_\b)-\phi\psi(v_\a)[\psi(v_\b), \phi\psi^{-1}(h')]\\
&=&[\psi\phi^{-1}(h'), \phi(v_\a)]\phi\psi(v_\b)+\phi\psi(v_\a)[h', \phi(v_\b)]\\
&=&[\phi^{-1}(h), \phi(v_\a)]\phi\psi(v_\b)+\phi\psi(v_\a)[\psi^{-1}(h), \phi(v_\b)]\\
&=&\a\phi^{-1}(h)\phi\psi(v_\a)\phi\psi(v_\b)+\b\psi^{-1}(h)\phi\psi(v_\a)\phi\psi(v_\b)\\
&=&(\a\phi^{-1}+\b\psi^{-1})(h)\phi\psi(v_\a)\phi\psi(v_\b).
\end{eqnarray*}
That is $A_{\a}A_{\b}\subset A_{\a\phi^{-1}+\b\psi^{-1}}$.
  \hfill $\square$

By Lemma 2.6, we can assert that
\begin{eqnarray*}
 A_{\a, \overline{i}}A_{\b, \overline{j}}\subset A_{\a\phi^{-1}+\b\psi^{-1}, \overline{i}+\overline{j}},~~\overline{i},\overline{j}\in \mathbb{Z}_2.
\end{eqnarray*}
In what follows, $A$ denotes a split regular BiHom-Poisson superalgebra and
\begin{eqnarray*}
A=H\oplus (\bigoplus_{\a\in \Lambda}A_{\a})
\end{eqnarray*}
the corresponding root spaces decomposition. Given a linear functional $\a:H_{\overline{0}}\rightarrow \mathbb{K}$, we denote by $-\a:H_{\overline{0}}\rightarrow \mathbb{K} $ the element in $H^{\ast}_{\overline{0}}$ defined by $(-\a)(h_{\overline{0}}):=-\a(h_{\overline{0}})$. We write
\begin{eqnarray*}
-\Lambda:=\{-\a:\a\in \Lambda\} ~~~\mbox{and}~~ \pm \Lambda:\Lambda\cup (-\Lambda).
\end{eqnarray*}
\begin{example}
Let $A=H\oplus (\bigoplus_{\a\in \Gamma}A_{\a})$ be a split Possion superalgebra, $\phi, \psi: A\rightarrow A$ two automorphism such that $\phi(H)=\psi(H)=H$ and $\phi\circ\psi=\psi\circ\phi$. By Example  2.2, we know that $(A, \mu', [\c,\c]', \phi, \psi)$ is a regular BiHom-Possion superalgebra. Then we have
\begin{eqnarray*}
A=H\oplus (\bigoplus_{\a\in \Gamma}A_{\a\psi^{-1}})
\end{eqnarray*}
makes of the regular BiHom-Possion superalgebra $(A, \mu', [\c,\c]', \phi, \psi)$ being the roots system $\Lambda=\{\a\psi^{-1}:\a\in \Gamma\}$.
\end{example}
\begin{definition}
 Let $\a,\b\in \Lambda$. We will say that $\a$ is connected to $\b$ if either
\begin{eqnarray*}
\b=\epsilon\a\phi^{z_1}\psi^{z_2}~~\mbox{for some $z_1,z_2\in \mathbb{Z}$ and $\epsilon\in \{-1,1\}$}
\end{eqnarray*}
or there exists $\{\a_1,\cdot\cdot\cdot,\a_k\}\subset \pm\Lambda$ with $k\geq 2$, such that \\
1. $\a_{1}\in \{\a \phi^{-n}\psi^{-r}:~~n,r \in \mathbb{N}\}$.\\
2. $\a_1\phi^{-1}+\a_2\psi^{-1}\in \pm\Lambda$,

 $\a_1\phi^{-2}+\a_2\phi^{-1}\psi^{-1}+\a_3\psi^{-1}\in \pm\Lambda$,

 $\a_1\phi^{-3}+\a_2\phi^{-2}\psi^{-1}+\a_3\phi^{-1}\psi^{-1}+\a_4\psi^{-1}\in \pm\Lambda$,

 $\cdot\cdot\cdot\cdot\cdot$

 $\a_1\phi^{-i}+\a_2\phi^{-i+1}\psi^{-1}+\a_3\phi^{-i+2}\psi^{-1}+\cdot\cdot\cdot+\a_i\phi^{-1}\psi^{-1}+\a_{i+1}\psi^{-1}\in \pm\Lambda$,

 $\a_1\phi^{-k+2}+\a_2\phi^{-k+3}\psi^{-1}+\a_3\phi^{-k+4}\psi^{-1}+\cdot\cdot\cdot+\a_{k-2}\phi^{-1}\psi^{-1}+\a_{k-1}\psi^{-1}\in \pm\Lambda$.\\
 3. $\a_1\phi^{-k+1}+\a_2\phi^{-k+2}\psi^{-1}+\a_3\phi^{-k+3}\psi^{-1}+\cdot\cdot\cdot+\a_{i}\phi^{-k+i}\psi^{-1}+\cdot\cdot\cdot+\a_{k-1}\phi^{-1}\psi^{-1}+\a_{k}\psi^{-1} \in \{\pm \b\phi^{-m}\psi^{-s}:m,s\in \mathbb{N}\}$.

 We will also say that  $\{\a_1,\cdot\cdot\cdot,\a_k\}$ is a connection from $\a$ to $\b$.
\end{definition}
The proof of the next result is analogous to the one of \cite{Zhang2018}. For the sake of completeness,  we give a sketch of the proof.

\begin{proposition}
The relation $\sim$ in $\Lambda$, defined by $\alpha\sim \b$ if and only if  $\a$ is connected to $\b$, is an equivalence relation.
\end{proposition}
 {\bf Proof.} If $\a\sim \b$, then either $\b=\epsilon\a\phi^{z_1}\psi^{z_2}$ for some $z_1,z_2\in \mathbb{Z}$ and $\epsilon\in \{-1,1\}$, and so $\b$
  is connected to $\a$; or there exists $\{\a_1,\cdot\cdot\cdot,\a_k\}\subset \pm\Lambda$ with $k\geq 2$, from $\a$ to $\b$ with
  \begin{eqnarray*}
\a_1\phi^{-k+1}+\a_2\phi^{-k+2}\psi^{-1}+\a_3\phi^{-k+3}\psi^{-1}+\cdot\cdot\cdot+\a_{i}\phi^{-k+i}\psi^{-1}+\cdot\cdot\cdot+\a_{k}\psi^{-1} =\epsilon\b\phi^{-m}\psi^{-s},
  \end{eqnarray*}
  for some $m,s\in \mathbb{N}$, $\epsilon\in \{-1,1\}$. Then we can verify that
 \begin{eqnarray*}
\{\b\phi^{-m}\psi^{-s}, -\epsilon\a_k\phi^{-1}, -\epsilon\a_{k-1}\phi^{-3}, -\epsilon\a_{k-2}\phi^{-5},\cdot\cdot\cdot,-\epsilon\a_{2}\phi^{-2k+3}\}
 \end{eqnarray*}
 is a connection from $\b$ to $\a$ and the relation $\sim$ is symmetric.

  Finally, suppose $\a\sim \b$ and $\b\sim \g$. If $\b=\epsilon\a\phi^{z_1}\psi^{z_2}$ for some $z_1,z_2\in \mathbb{Z}$, $\epsilon\in \{-1,1\}$ and
  $\g=\epsilon'\a\phi^{z_3}\psi^{z_4}$ for some $z_3,z_4\in \mathbb{Z}$, $\epsilon'\in \{-1,1\}$, it is clear that $\a\sim \g$.
  Hence  suppose $\{\a_1,\cdot\cdot\cdot,\a_k\}$ with $k\geq 2$ is a connection from $\a$ to $\b$ which satisfies
  \begin{eqnarray*}
 \a_1\phi^{-k+1}+\a_2\phi^{-k+2}\psi^{-1}+\cdot\cdot\cdot+\a_k\psi^{-1}=\epsilon\b\phi^{-m}\psi^{-s}
  \end{eqnarray*}
 for some $m,s\in \mathbb{N}$, $\epsilon\in \{-1,1\}$, and $\{h_1,\cdot\cdot\cdot,h_p\}$ is a connection from $\b$ to $\g$. Then
 $\{\a_1,\cdot\cdot\cdot,\a_k,  \epsilon h_2,\cdot\cdot\cdot, \epsilon h_p\}$   is connection from $\a$ to $\g$, so the connection relation is also transitive.
  \hfill $\square$

  By  Proposition 2.9 we can consider the quotient set
  \begin{eqnarray*}
 \Lambda/\sim=\{[\a]:\a\in \Lambda\},
  \end{eqnarray*}
with $[\a]$ being the set of nozero roots which are connected to $\a$.
Our next goal is to associate an ideal $I_{[\a]}$ to $[\a]$. Fix $[\a]\in \Lambda/\sim$, we start by defining
\begin{eqnarray*}
I_{H,[\a]}=span_{\mathbb{K}}\{[A_{\b\psi^{-1}}, A_{-\b\phi^{-1}}]+A_{\b\psi^{-1}}A_{-\b\phi^{-1}}:\b\in [\a]\}.
\end{eqnarray*}
Now we define
\begin{eqnarray*}
V_{[\a]}:=\bigoplus_{\b\in [\a]}A_{\b}.
\end{eqnarray*}
Finally, we denote by $I_[\a]$ the direct sum of the two subspaces above:
\begin{eqnarray*}
I_{[\a]}:=I_{H,[\a]}\oplus V_{[\a]}.
\end{eqnarray*}

\begin{proposition} For any $[\a]\in \Lambda/\sim$, the following assertions hold.

(1). $[I_{[\a]},I_{[\a]}]+I_{[\a]}I_{[\a]}\subset I_{[\a]}$,

(2). $\phi(I_{[\a]})=I_{[\a]}$ and $\psi(I_{[\a]})=I_{[\a]}$,

(3). For any $[\b]\neq [\a]$, we have $[I_{[\a]}, I_{[\b]}]+I_{[\a]}I_{[\b]}=0$.
\end{proposition}
{\bf Proof.} (1) First we check that $[I_{[\a]},I_{[\a]}]\subset I_{[\a]}$, we can write
\begin{eqnarray}
 [I_{[\a]},I_{[\a]}]&=&[I_{H,[\a]}\oplus V_{[\a]},I_{H,[\a]}\oplus V_{[\a]}]\nonumber\\
 &\subset&[I_{H,[\a]}, V_{[\a]}]+[V_{[\a]}, I_{H,[\a]}]+[V_{[\a]}, V_{[\a]}],
\end{eqnarray}
Given $\b\in [\a]$, we have $[I_{H,[\a]}]\subset A_{\b\psi^{-1}}$. Since $\b\psi^{-1}\in [\a]$, it follows that $[I_{H,[\a]}, A_{\b}]\subset V_{[\a]}$.

By a similar argument, we get $[A_{\b}, I_{H,\a}]\subset V_{[\a]}$.

Next we consider  $[V_{[\a]},V_{[\a]}]$. If we take $\b,\g\in [\a]$ such that $[A_{\b}, A_{\g}]\neq 0$, then $[A_{\b}, A_{\g}]\subset A_{\b\phi^{-1}+\g\psi^{-1}}$.
If $\b\phi^{-1}+\g\psi^{-1}=0$, we get $[A_{\b}, A_{-\g}]\subset H$ and so $[A_{\b}, A_{-\g}]\subset I_{H,[\a]}$. Suppose that $\b\phi^{-1}+\g\psi^{-1}\in \Lambda$. We infer that $\{\b,\g\}$ is connection from $\b$ to $\b\phi^{-1}+\g\psi^{-1}$. The transitivity of $\sim$ now gives that $\b\phi^{-1}+\g\psi^{-1}\in [\a]$ and so $[A_{\b}, A_{\g}]\subset V_{[\a]}$. Hence
\begin{eqnarray}
[V_{[\a]},V_{[\a]}]\in I_{[\a]}.
\end{eqnarray}
 From (2.1) and (2.2), we get $[I_{[\a]},I_{[\a]}]\subset I_{[\a]}$.

Second, we will check that $I_{[\a]}I_{[\a]}\subset I_{[\a]}$. We have
\begin{eqnarray}
I_{[\a]}I_{[\a]}&=&(I_{H,[\a]}\oplus V_{[\a]})(I_{H,[\a]}\oplus V_{[\a]})\nonumber\\
&\subset& I_{H,[\a]}I_{H,[\a]}+I_{H,[\a]}V_{[\a]}+V_{[\a]}I_{H,[\a]}+V_{[\a]}V_{[\a]}.
\end{eqnarray}
Similar considerations, we have
\begin{eqnarray*}
I_{H,[\a]}V_{[\a]}+V_{[\a]}I_{H,[\a]}+V_{[\a]}V_{[\a]}\subset I_{H,[\a]}.
\end{eqnarray*}
Hence, it just remains to check that $I_{H,[\a]}I_{H,[\a]}$, observe that
\begin{eqnarray}
I_{H,[\a]}I_{H,[\a]}&\subset&(\sum_{\b\in [\a]}[A_{\b\psi^{-1}}, A_{-\b\phi^{-1}}]+A_{\b\psi^{-1}}A_{-\b\phi^{-1}})H\nonumber\\
&\subset& (\sum_{\b\in [\a]}[A_{\b\psi^{-1}}, A_{-\b\phi^{-1}}])H+(\sum_{\b\in [\a]}A_{\b\psi^{-1}}A_{-\b\phi^{-1}})H.
\end{eqnarray}
Consider the first summand on the right hand side of (2.4). By BiHom-Leibniz superidentity, we have
\begin{eqnarray*}
&&[A_{\b\psi^{-1}}, A_{-\b\phi^{-1}}]\phi\phi^{-1}(H)\\
&\subset& [A_{-\b\psi^{-1}}\phi^{-1}(H), \phi(A_{-\b\phi^{-1}})]+\phi(A_{\b\psi^{-1}})[\phi^{-1}(H), \phi\psi^{-1}(A_{-\b\phi^{-1}})]\\
&\subset& [A_{-\b\psi^{-1}\phi^{-1}}, A_{-\b\phi^{-1}\phi^{-1}}]+A_{\b\psi^{-1}\phi^{-1}})A_{-\b\phi^{-1}\psi\phi^{-2}}\\
&\subset& I_{H,[\a]}.
\end{eqnarray*}
Next we consider the last summand on the right hand side of (2.4). By BiHom-associativity, we know
\begin{eqnarray*}
(A_{\b\psi^{-1}}A_{-\b\phi^{-1}})\psi(\psi^{-1}(H))&=&\phi(A_{\b\psi^{-1}})(A_{-\b\phi^{-1}}\psi^{-1}(H))\\
&\subset& A_{\b\psi^{-1}\phi^{-1}}A_{-\b\phi^{-1}\psi^{-1}}\subset I_{H,[\a]}.
\end{eqnarray*}

(2) It is easy to check that $\phi(I_{[\a]})=I_{[\a]}$ and $\psi(I_{[\a]})=I_{[\a]}$.

(3) We will study the expression $[I_{[\a]},I_{[\b]}]+I_{[\a]}I_{[\b]}$. Observe that
\begin{eqnarray}
[I_{[\a]},I_{[\b]}]&=&[I_{H,[\a]}\oplus V_{[\a]},I_{H,[\b]}\oplus V_{[\b]}]\nonumber\\
 &\subset&[I_{H,[\a]}, V_{[\b]}]+[ V_{[\a]},I_{H,[\b]} ]+[V_{[\a]}, V_{[\b]}],
\end{eqnarray}
and
\begin{eqnarray}
I_{[\a]}I_{[\b]}&=&(I_{H,[\a]}\oplus V_{[\a]})(I_{H,[\b]}\oplus V_{[\b]})\nonumber\\
&\subset&  I_{H,[\a]}I_{H,[\b]}+I_{H,[\a]}V_{[\b]}+V_{[\a]}I_{H,[\b]}+V_{[\a]}V_{[\b]}.
\end{eqnarray}

First we  consider $[V_{[\a]}, V_{[\b]}]+V_{[\a]}V_{[\b]}$ and suppose that there exist $\a_1\in [\a]$ and  $\b_1\in [\b]$ such that $[A_{\a_1}, A_{\b_1}]+A_{\a_1}A_{\b_1}\neq 0$. As necessarily $\a_1\phi^{-1}\neq-\b_1\psi^{-1}$, then $\a_1\phi^{-1}+\b_1\psi^{-1}\in \Lambda$. So $\{\a_1,\b_1, -\a_1\phi^{-1}\}$ is a connection between $\a_1$ and $\b_1$. By the transitivity of the connection relation we see $\a\in [\b]$, a contradicition. Hence $[A_{\a_1}, A_{\b_1}]+A_{\a_1}A_{\b_1}=0$, and so
\begin{eqnarray*}
[V_{[\a]}, V_{[\b]}]+V_{[\a]}V_{[\b]}=0.
\end{eqnarray*}
Next we consider the first summand $[I_{H,[\a]}, V_{[\b]}]$ on the right hand side of (2.5) and the second one $I_{H,[\a]}V_{[\b]}$ of (2.6), and suppose that there exist $\a_1\in [\a]$ and $\b_1\in [\b]$ such that
\begin{eqnarray*}
[[A_{\a_1}, A_{-\a_1}], A_{\b_1}]+[A_{\a_1}A_{-\a_1}, A_{\b_1}]+[A_{\a_1}, A_{-\a_1}] A_{\b_1}+(A_{\a_1}A_{-\a_1}) A_{\b_1}\neq 0.
\end{eqnarray*}
Then some of the four sunmands are different from zero.

If
\begin{eqnarray*}
[[A_{\a_1}, A_{-\a_1}], A_{\b_1}]\neq 0,
\end{eqnarray*}
then  BiHom-Leibniz identity gives
\begin{eqnarray*}
0&\neq& [[A_{\a_1}, A_{-\a_1}], \phi\phi^{-1}(A_{\b_1})]\\
&\subset& [[A_{\a_1}, \phi^{-1}(A_{\b_1})], \phi(A_{-\a_1})]+[\phi(A_{\a_1}), [A_{-\a_1}, \psi^{-1}(A_{\b_1})]]\\
&\subset&[[A_{\a_1}, \phi^{-1}(A_{\b_1})], \phi(A_{-\a_1})]+[[\phi^{-1}\psi(A_{-\a_1}), \phi^{-1}(A_{\b_1})], \phi^2\psi^{-1}(A_{\a_1})].
\end{eqnarray*}
Hence
\begin{eqnarray*}
[A_{\a_1}, \phi^{-1}(A_{\b_1})]+[\phi^{-1}\psi(A_{-\a_1}), \phi^{-1}(A_{\b_1})]\neq 0,
\end{eqnarray*}
which contradicts (2.6). Therefore, $[[A_{\a_1}, A_{-\a_1}], A_{\b_1}]=0$.

If the second, third or fouth summand were nonzero, we can argue as above but using the BiHom-Leibniz or BiHom-associativity superidentites to show that these products are zero. Consequently,
\begin{eqnarray*}
[I_{H,[\a]},V_{[\b]}] +I_{H,[\a]}V_{[\b]}=0.
\end{eqnarray*}
In a similar way we can prove that the remaining summands in (2.5) and (2.6) are zero, and the proof is complete. \hfill $\square$

\begin{proposition} For any $[\a]\in \Lambda/\sim$, we have
\begin{eqnarray*}
I_{H,[\a]}H+HI_{H[\a]}\in I_{H[\a]}.
\end{eqnarray*}
\end{proposition}
{\bf Proof.}  Fix any $\b\in [\a]$. On the one hand, by the BiHom-Leibniz superidentity, we get
\begin{eqnarray*}
[A_\b, A_{-\b}]H+H[A_\b, A_{-\b}]\in I_{H[\a]}.
\end{eqnarray*}
And on the other hand, by BiHom-associativity we know
\begin{eqnarray*}
(A_\b A_{-\b})H+H(A_\b A_{-\b})\in I_{H[\a]}.
\end{eqnarray*}
\begin{theorem}
(1) For any $[\a]\in \Lambda/\sim$, the linear space $I_{[\a]}=I_{H,[\a]}+V_{[\a]}$ of $A$ associated to $[\a]$ is an ideal of $A$.

(2) If $A$ is simple, then there exists a connection from $\a$ to $\b$ for any $\a,\b\in \Lambda$ and $H=\sum_{\a\in \Lambda}([A_{\a\psi^{-1}}, A_{-\a\phi^{-1}}]+A_{\a\psi^{-1}}A_{-\a\phi^{-1}})$.
\end{theorem}
{\bf Proof.} (1) Since $[I_{[\a]},H]\subset I_{[\a]}$, by  (2.4) and  (2.5) we have
\begin{eqnarray*}
[I_{[\a]}, A]=[I_{[\a]}, H\oplus(\bigoplus_{\b\in [\a]} A_{\b})\oplus (\bigoplus_{\g\in [\a]} A_{\g})]\subset I_{[\a]}.
\end{eqnarray*}
According to Proposition 2.9 and Proposition 2.10,  we have
\begin{eqnarray*}
&&I_{[\a]} A+AI_{[\a]}\\
&=&I_{[\a]}(H\oplus(\bigoplus_{\b\in [\a]} A_{\b})\oplus (\bigoplus_{\g\notin [\a]} A_{\g}))\\
&&+(H\oplus(\bigoplus_{\b\in [\a]} A_{\b})\oplus (\bigoplus_{\g\notin [\a]} A_{\g}))I_{[\a]}\subset I_{[\a]}.
\end{eqnarray*}
As we also have $\phi(I_{[\a]})=I_{[\a]}$ and $\psi(I_{[\a]})=I_{[\a]}$. So we conclude that $I_{[\a]}$ is an ideal of $A$.

(2) The simplicity of $A$ implies $I_{[\a]}=A$. From here, it is clear that $[\a]=\Lambda$ and $H=\sum_{\a\in \Lambda}([A_{\a\psi^{-1}}, A_{-\a\phi^{-1}}]+A_{\a\psi^{-1}}A_{-\a\phi^{-1}})$.
\hfill $\square$

\begin{theorem} We have
\begin{eqnarray*}
A=U+\sum_{[\a]\in \Lambda/\sim}I_{[\a]},
\end{eqnarray*}
where $U$ is a linear complement in $H$ of $span_{\mathbb{K}}\{[A_{\a\psi^{-1}}, A_{-\a\phi^{-1}}]+A_{\a\psi^{-1}}A_{-\a\phi^{-1}} :\a\in \Lambda \}$ and any $I_{[\a]}$ is one of the ideals of $A$ described in Theorem 2.12, satisfying $[I_{[\a]},I_{[\b]}]+ I_{[\a]}I_{[\b]}=0$ if $[\a]\neq[\b]$.
\end{theorem}
{\bf Proof.}  $I_{[\a]}$ is well defined and is an ideal of $A$, being clear that
 \begin{eqnarray*}
 A=H\oplus\sum_{[\a]\in \Lambda} A_{[\a]}=U+\sum_{[\a]\in \Lambda/\sim}I_{[\a]}.
 \end{eqnarray*}
 Finally, Proposition 2.10 gives us $[I_{[\a]},I_{[\b]}]+I_{[\a]}I_{[\b]}=0$ if $[\a]\neq[\b]$.\hfill $\square$

 Let us denote by $Z(A):=\{v\in A:[v,A]+vA+Av=0\}$ the center of $A$.

 \begin{corollary}  If $Z(A)=0$ and $H=\sum_{\a\in \Lambda}([A_{\a\psi^{-1}}, A_{-\a\phi^{-1}}]+A_{\a\psi^{-1}}A_{-\a\phi^{-1}})$. Then $A$ is the direct sum of the ideals given in Theorem 2.12,
 \begin{eqnarray*}
A=\bigoplus_{[\a]\in \Lambda/\sim}I_{[\a]},
\end{eqnarray*}
Furthermore, $[I_{[\a]},I_{[\b]}]+I_{[\a]}I_{[\b]}=0$ if $[\a]\neq[\b]$.\hfill $\square$.
\end{corollary}
{\bf Proof.} Since $H=\sum_{\a\in \Lambda}([A_{\a\psi^{-1}}, A_{-\a\phi^{-1}}]+A_{\a\psi^{-1}}A_{-\a\phi^{-1}})$, it follows that  $A=\sum_{[\a]\in \Lambda/\sim}I_{[\a]}$. To verify the direct character of the sum, take some $v\in I_{[\a]}\cap(\sum_{[\b]\in\Lambda/\sim,[\b]\neq[\a]}I_{[\b]})$. Since $v\in I_{[\a]}$, the fact $[I_{[\a]},I_{[\b]}]+I_{[\a]}I_{[\b]}=0$ when  $[\a]\neq[\b]$ gives us
\begin{eqnarray*}
[v, \sum_{[\b]\in\Lambda/\sim,[\b]\neq[\a]}I_{[\b]}]+v(\sum_{[\b]\in\Lambda/\sim,[\b]\neq[\a]}I_{[\b]})+(\sum_{[\b]\in\Lambda/\sim,[\b]\neq[\a]}I_{[\b]})v=0.
\end{eqnarray*}
  In a similar way,  $v\in \sum_{[\b]\in\Lambda/\sim,[\b]\neq[\a]}I_{[\b]}$ implies $[v,I_{[\a]}]+v I_{[\a]}+I_{[\a]}v=0$. That is $v\in Z(A)$ and so $v=0$.              \hfill $\square$
\section{The simple components}
\def\theequation{\arabic{section}. \arabic{equation}}
\setcounter{equation} {0}

In this section we focus on the simplicity of split regular BiHom-Poisson superalgebras by
centering our attention on those of maximal length, we recall that a roots system $\Lambda$ of a split regular BiHom-Poisson superalgebra $A$ is called symmetric if it is satisfies that $\a\in \Lambda$ implies $-\a\in \Lambda$. From now on we will suppose that $\Lambda$ is symmetric.

For the grading of $I$, we have
\begin{eqnarray*}
I=I_{\overline{0}}\oplus I_{\overline{1}}=((I_{\overline{0}}\cap H_{\overline{0}})\oplus (\oplus_{\a\in \Lambda} (I_{\overline{0}}\cap L_{\a, \overline{0}})))\oplus ((I_{\overline{1}}\cap H_{\overline{1}})\oplus (\oplus_{\a\in \Lambda} (I_{\overline{1}}\cap L_{\a, \overline{1}}))).
\end{eqnarray*}

\begin{lemma}
 Suppose $H=\sum_{\a\in \Lambda}([A_{\a\psi^{-1}}, A_{-\a\phi^{-1}}]+A_{\a\psi^{-1}}A_{-\a\phi^{-1}})$. If $I$ is an ideal of $A$ such that $I\subset H$, then $I\subset Z(A)$.
\end{lemma}
{\bf Proof.} Observe that $[I,H]\subset [H,H]=0$ and
\begin{eqnarray*}
&&[I, \bigoplus_{\a\in \Lambda}A_{\a}] +I(\bigoplus_{\a\in \Lambda}A_{\a})+(\bigoplus_{\a\in \Lambda}A_{\a})I\\
&\subset& I\cap (\bigoplus_{\a\in \Lambda}A_{\a})\subset H\cap (\bigoplus_{\a\in \Lambda}A_{\a})=0.
\end{eqnarray*}
Since $H=\sum_{\a\in \Lambda}([A_{\a\psi^{-1}}, A_{-\a\phi^{-1}}]+A_{\a\psi^{-1}}A_{-\a\phi^{-1}})$, by the BiHom-Leibniz superidentity and the above observation, we obtain $HI+IH=0$. So $I\subset Z(A)$.
\hfill $\square$

 We  use the same notions  of  \cite{Calderon2016} and \cite{Zhang2018},   denote by $\Lambda_{\overline{0}}:=\{\a \in \Lambda: L_{\a, \overline{0}}\neq 0\}$ and $\Lambda_{\overline{1}}:=\{\a \in \Lambda: L_{\a, \overline{1}}\neq 0\}$, so $\Lambda=\Lambda_{\overline{0}}\cup \Lambda_{\overline{1}}$.

\begin{definition} A split regular BiHom-Poisson superalgebra $A$ is root multiplicative if  $\a\in\Lambda_{\overline{i}} ,\b\in \Lambda_{\overline{i}}, \overline{i},\overline{j}\in \mathbb{Z}_2 $  such that $\a\psi^{-1}+\b\phi^{-1}\in \Lambda_{\overline{i}+\overline{j}}$, then $[A_{\a, \overline{i}}, A_{\b, \overline{j}}]+A_{\a, \overline{i}}A_{\b, \overline{j}}\neq 0$.
\end{definition}

\begin{definition}
A split regular BiHom-Poisson superalgebra $A$ is of maximal length if dim$A_{\a, \overline{i}}=1$ for any $\a\in \Lambda$ and $\overline{i}\in \mathbb{Z}_2$.
\end{definition}

Observe that if $A$ is of  maximal lenth, then we have
\begin{eqnarray}
I=((I_{\overline{0}}\cap H_{\overline{0}})\oplus (\oplus_{\a\in \Lambda'_{\overline{0}}}  L_{\a, \overline{0}}))\oplus ((I_{\overline{1}}\cap H_{\overline{1}})\oplus (\oplus_{\a\in \Lambda'_{\overline{1}}}  L_{\a, \overline{1}})),
\end{eqnarray}
where $\Lambda'_{\overline{i}}=\{\a\in \Lambda_{\overline{i}}: I_{\overline{i}}\cap L_{\a, \overline{i}}\neq 0\}, \overline{i}\in \mathbb{Z}_2$.

\begin{theorem} Let $A$ be a  split regular BiHom-Poisson superalgebra  of maximal length and root  multiplicative. Then $A$ is simple if and only if $Z(A)=0$, $H=\sum_{\a\in \Lambda}([A_{\a\psi^{-1}}, A_{-\a\phi^{-1}}]+A_{\a\psi^{-1}}A_{-\a\phi^{-1}})$ and $\Lambda$ has all of its elements connected.
\end{theorem}
{\bf Proof.} Suppose $A$ is simple. Since $Z(A)$ is an ideal of $A$, we have $Z(A)=0$. Now Theorem 2.12(2) completes the proof of the direct implication.

To prove the converse, consider a nonzero ideal of $A$. By (3.1), we can write $I=((I_{\overline{0}}\cap H_{\overline{0}})\oplus (\oplus_{\a\in \Lambda'_{\overline{0}}}  L_{\a, \overline{0}}))\oplus ((I_{\overline{1}}\cap H_{\overline{1}})\oplus (\oplus_{\a\in \Lambda'_{\overline{1}}}  L_{\a, \overline{1}}))$, where $\Lambda'_{\overline{i}}\subset \Lambda_{\overline{i}}, \overline{i}\in \mathbb{Z}_2$, and some  $\Lambda'_{\overline{i}}\neq \emptyset$ as consequence of Lemma 3.1. Let us fix some $\a_0\in \Lambda'_{\overline{i}}$ with $0\neq A_{\a_0, \overline{i}}\subset I$. Since $\phi(I)=I$ and $\psi(I)=I$, it follows  that
\begin{eqnarray}
\mbox{if $\a\in \Lambda'_{\overline{i}}$, then}~~ \{\a\phi^{z_1}\psi^{z_2}:z_1,z_2\in \mathbb{Z}\}\subset \Lambda'_{\overline{i}}.
\end{eqnarray}
In particular,
\begin{eqnarray}
\{A_{\a_0\phi^{z_1}\psi^{z_2}, \overline{i}}: z_1,z_2\in \mathbb{Z}\}\subset I.
\end{eqnarray}

Now, let us take any $\b\in \Lambda$ satisfying $\b\notin \{\a_0\phi^{z_1}\psi^{z_2}:z_1,z_2\in \mathbb{Z}\}$. Since $\a_0$ and $\b$ are connected, we have a connection $\{\a_1,\cdot\cdot\cdot,\a_k\}, k\geq 2$, from $\a_0$ to $\b$ satisfying:

 $\a_{1}=\a_0 \phi^{-n}\psi^{-r}:~~n,r \in \mathbb{N}$.

 $\a_1\phi^{-1}+\a_2\psi^{-1}\in \Lambda$,

 $\a_1\phi^{-2}+\a_2\phi^{-1}\psi^{-1}+\a_3\psi^{-1}\in \Lambda$,

 $\a_1\phi^{-3}+\a_2\phi^{-2}\psi^{-1}+\a_3\phi^{-1}\psi^{-1}+\a_4\psi^{-1}\in \Lambda$,

 $\cdot\cdot\cdot\cdot\cdot$

 $\a_1\phi^{-i}+\a_2\phi^{-i+1}\psi^{-1}+\a_3\phi^{-i+2}\psi^{-1}+\cdot\cdot\cdot+\a_i\phi^{-1}\psi^{-1}+\a_{i+1}\psi^{-1}\in \Lambda$,

 $\a_1\phi^{-k+2}+\a_2\phi^{-k+3}\psi^{-1}+\a_3\phi^{-k+4}\psi^{-1}+\cdot\cdot\cdot+\a_{k-2}\phi^{-1}\psi^{-1}+\a_{k-1}\psi^{-1}\in \Lambda$,

 $\a_1\phi^{-k+1}+\a_2\phi^{-k+2}\psi^{-1}+\a_3\phi^{-k+3}\psi^{-1}+\cdot\cdot\cdot+\a_{i}\phi^{-k+i}\psi^{-1}+\cdot\cdot\cdot+\a_{k-1}\phi^{-1}\psi^{-1}+\a_{k}\psi^{-1} =\epsilon\b\phi^{-m}\psi^{-s}:m,s\in \mathbb{N}$.

Taking into account that $\a_1,\a_2\in \Lambda$ and $\a_1\phi^{-1}+\a_2\psi^{-1}\in \Lambda$, we have $\a_1\in \Lambda_{\overline{i}}$ and $\a_1\in \Lambda_{\overline{j}}$ such that $\a_1\phi^{-1}+\a_2\psi^{-1}\in \Lambda_{\overline{i}+\overline{j}}$.   The root multiplicativity and maximal length of $A$ allow us to assert that either $0\neq [A_{\a_1, \overline{i}}, A_{\a_2, \overline{j}}]=A_{\a_1\phi^{-1}+\a_2\psi^{-1}, \overline{i}+\overline{j}}$ or $0\neq A_{\a_1, \overline{i}}A_{\a_2, \overline{j}}+A_{\a_2, \overline{j}}A_{\a_1, \overline{i}}=A_{\a_1\phi^{-1}+\a_2\psi^{-1}, \overline{i}+\overline{j}} $.

Since $0\neq A_{\a_1, \overline{i}}\subset I$ as a consequence of (3.3) we get
\begin{eqnarray*}
0\neq A_{\a_1\phi^{-1}+\a_2\psi^{-1}, \overline{i}+\overline{j}}\subset I.
\end{eqnarray*}
A similar argument applied to $\a_1\phi^{-1}+\a_2\psi^{-1}, \a_3$, and
\begin{eqnarray*}
(\a_1\phi^{-1}+\a_2\psi^{-1})\phi^{-1} +\a_3\psi^{-1}=\a_1\phi^{-2}+\a_2\phi^{-1}\psi^{-1} +\a_3\psi^{-1}
\end{eqnarray*}
gives us $0\neq A_{\a_1\phi^{-2}+\a_2\phi^{-1}\psi^{-1} +\a_3\psi^{-1}, \overline{k}}\subset I$. We can follow this process with the connection $\{\a_1,\cdot\cdot\cdot,\a_k\}$ to get
\begin{eqnarray*}
0\neq A_{\a_1\phi^{-k+1}+\a_2\phi^{-k+2}\psi^{-1}+\cdot\cdot\cdot+\a_{k}\psi^{-1}, \overline{h} }\subset I,
\end{eqnarray*}
and then
\begin{eqnarray*}
\mbox{either}~~ A_{\b\phi^{-m}\psi^{-s}, \overline{h}}\subset I~~ \mbox{or}~~ A_{-\b\phi^{-m}\psi^{-s}, \overline{h}}\subset I.
\end{eqnarray*}
From (3.2) and (3.3), we have
\begin{eqnarray*}
\mbox{either}~~ A_{\a\phi^{-z_1}\psi^{-z_2}, \overline{i}}\subset I ~~\mbox{or}~~ A_{-\a\phi^{-z_1}\psi^{-z_2}, \overline{i}}\subset I.
\end{eqnarray*}
This can be reformulated by saying that for any $\a\in \Lambda$,  either $\{\a\phi^{-z_1}\psi^{-z_2}\}$ or $\{-\a\phi^{-z_1}\psi^{-z_2}\}$ is contained in $\Lambda_{I}$. Taking now into account $H=\sum_{\a\in \Lambda}([A_{\a\psi^{-1}}, A_{-\a\phi^{-1}}]+A_{\a\psi^{-1}}A_{-\a\phi^{-1}})$, we have
\begin{eqnarray}
H\subset I.
\end{eqnarray}
Now for any $\a\in \Lambda$, since $A_{\a}=[H, A_{\a\psi}]$ by the maximal length of $A$, (3.4) gives us $A_{\a}\subset I$ and so $A=I$. That is, $A$ is simple. \hfill $\square$

\begin{theorem}Let $A$ be a  split regular BiHom-Poisson superalgebra  of maximal length and root  multiplicative with $Z(A)=0$ and satisfying $H=\sum_{\a\in \Lambda}([A_{\a\psi^{-1}}, A_{-\a\phi^{-1}}]+A_{\a\psi^{-1}}A_{-\a\phi^{-1}})$. Then
\begin{eqnarray*}
A=\bigoplus_{[\a]\in \Lambda/\sim}I_{[\a]},
\end{eqnarray*}
where any $I_{[\a]}$ is a simple split ideal having its roots system $\Lambda_{I_{[\a]}}$, with all of its elements  $\Lambda_{I_{[\a]}}$-connected.
\end{theorem}
{\bf Proof.} By Corollary 2.11, we can write $A$ as the direct sum $\bigoplus_{[\a]\in \Lambda/\sim}I_{[\a]}$ of the family of ideals
\begin{eqnarray*}
I_{[\a]}=I_{H,[\a]}\oplus V_{[\a]}=span_{\mathbb{K}}\{[A_{\b\psi^{-1}}, A_{-\b\phi^{-1}}]+A_{\b\psi^{-1}}A_{-\b\phi^{-1}}:\b\in [\a]\}\oplus_{\b\in [\a]} A_{\b},
\end{eqnarray*}
where each $I_{[\a]}$ is a split regular BiHom-Poisson superalgebra  with root system $A_{I_{[\a]}}=[\a]$. To make use of Theorem 3.4 in each $I_{[\a]}$, we observe that the root multiplicativity of $A$ and Proposition 2.10 show that $A_{I_{[\a]}}$ has all of its elements  $A_{I_{[\a]}}$ connected, that is, connected through connections contained in $A_{I_{[\a]}}$. Moreover, each $I_{[\a]}$ is root multiplicative by the root multiplicativity of $A$. So we infer that  $I_{[\a]}$ is of maximal length, and finally its center $Z(I_{[\a]})=0$. As consequence $[I_\a,I_\b]+I_\a I_\b=0$ if $[\a]\neq[\b]$. Applying Theorem 3.4, we conclude  that $I_{[\a]}$ is simple and $A=\bigoplus_{[\a]\in \Lambda/\sim}I_{[\a]}$. \hfill $\square$

 \begin{center}
 {\bf ACKNOWLEDGEMENT}
 \end{center}

  The paper is supported by  the NSF of China (No. 11761017) and the Youth Project for Natural Science Foundation of Guizhou provincial department of education (No. KY[2018]155).

\end{document}